\documentclass[12pt,leqno]{article}
\usepackage{amsfonts}
\usepackage{xcolor}
\usepackage{amsmath, amsthm, amsfonts, amssymb, color}
\usepackage{mathrsfs}
\usepackage{color}
\usepackage[notcite,notref]{showkeys}

\theoremstyle{definition}

\newcommand{\scr}[1]{\mathscr #1}
\definecolor{wco}{rgb}{0.5,0.2,0.3}

\numberwithin{equation}{section} \theoremstyle{remark}

\newcommand{\ua}{\uparrow}

\title{{\bf  Singular Degenerate SDEs: Well-Posedness and Exponential Ergodicity }\footnote{Feng-Yu Wang is supported in part by the National Key R\&D Program of China (2022YFA1006000, 2020YFA0712900) and NNSFC (11831014, 11921001).  Martin Grothaus and Panpan Ren acknowledge support by the Alexander von Humboldt Foundation.}}
\author{
{\bf Panpan Ren$^{c)}$,   Martin Grothaus$^{b)}$,    Feng-Yu Wang$^{a)}$  }\\
\footnotesize{$^{a)}$ Center for Applied Mathematics, Tianjin University, Tianjin 300072, China}\\
\footnotesize{$^{b)}$ Mathematics Department,  University of Kaiserslautern-Landau,  Kaiserslautern  67653, Germany}\\
\footnotesize{$^{c)}$ Mathematics Department,  City University of Hong Kong,  Hong Kong, Kowloon,  China. }\\
\footnotesize{ grothaus@mathematik.uni-kl.de,  panparen@cityu.edu.hk,  wangfy@tju.edu.cn}}
\begin{document}
\allowdisplaybreaks
\def\R{\mathbb R}  \def\ff{\frac} \def\ss{\sqrt} \def\B{\mathbf
B}
\def\N{\mathbb N} \def\kk{\kappa} \def\m{{\bf m}}
\def\ee{\varepsilon}\def\ddd{D^*}
\def\dd{\delta} \def\DD{\Delta} \def\vv{\varepsilon} \def\rr{\rho}
\def\<{\langle} \def\>{\rangle}
  \def\nn{\nabla} \def\pp{\partial} \def\E{\mathbb E}
\def\d{\text{\rm{d}}} \def\bb{\beta} \def\aa{\alpha} \def\D{\scr D}
  \def\si{\sigma} \def\ess{\text{\rm{ess}}}\def\s{{\bf s}}
\def\beg{\begin} \def\beq{\begin{equation}}  \def\F{\scr F}
\def\Ric{\mathcal Ric} \def\Hess{\text{\rm{Hess}}}
\def\e{\text{\rm{e}}} \def\ua{\underline a} \def\OO{\Omega}  \def\oo{\omega}
 \def\tt{\tilde}\def\[{\lfloor} \def\]{\rfloor}
\def\cut{\text{\rm{cut}}} \def\P{\mathbb P} \def\ifn{I_n(f^{\bigotimes n})}
\def\C{\scr C}      \def\aaa{\mathbf{r}}     \def\r{r}
\def\gap{\text{\rm{gap}}} \def\prr{\pi_{{\bf m},\varrho}}  \def\r{\mathbf r}
\def\Z{\mathbb Z} \def\vrr{\varrho} \def\ll{\lambda}
\def\L{\scr L}\def\Tt{\tt} \def\TT{\tt}\def\II{\mathbb I}
\def\i{{\rm in}}\def\Sect{{\rm Sect}}  \def\H{\mathbb H}
\def\M{\mathbb M}\def\Q{\mathbb Q} \def\texto{\text{o}} \def\LL{\Lambda}
\def\Rank{{\rm Rank}} \def\B{\scr B} \def\i{{\rm i}} \def\HR{\hat{\R}^d}
\def\to{\rightarrow}\def\l{\ell}\def\iint{\int}\def\gg{\gamma}
\def\EE{\scr E} \def\W{\mathbb W}
\def\A{\scr A} \def\Lip{{\rm Lip}}\def\S{\mathbb S}
\def\BB{\scr B}\def\Ent{{\rm Ent}} \def\i{{\rm i}}\def\itparallel{{\it\parallel}}
\def\g{{\mathbf g}}\def\Sect{{\mathcal Sec}}\def\T{\mathcal T}\def\BB{{\bf B}}
\def\f{\mathbf f} \def\g{\mathbf g}\def\BL{{\bf L}}  \def\BG{{\mathbb G}}
\def\Bd{{D^E}} \def\BdP{D^E_\phi} \def\Bdd{{\bf \dd}} \def\Bs{{\bf s}} \def\GA{\scr A}
\def\Bg{{\bf g}}  \def\Bdd{\psi_B} \def\supp{{\rm supp}}\def\div{{\rm div}}
\def\ddiv{{\rm div}}\def\osc{{\bf osc}}\def\1{{\bf 1}}\def\BD{\mathbb D}\def\GG{\Gamma}
\def\H{{\bf H}} \def\n{{\bf n}}
\maketitle

\begin{abstract} The well-posedness and exponential ergodicity are proved for stochastic Hamiltonian systems
containing a singular drift term which  is locally integrable in the component with noise. As an application, the well-posedness and uniform exponential ergodicity are derived for a class of singular degenerated McKean-Vlasov SDEs. \end{abstract} \noindent

 \vskip 2cm

 \section{Introduction}

Let $d_1,d_2\in \mathbb N$. For fixed $T\in (0,\infty]$, consider the following degenerate SDE for $(X_t,Y_t)\in\R^{d_1+d_2}=\R^{d_1}\times \R^{d_2}$:
\beq\label{E1}\beg{cases} \d X_t= Z_t^{(1)}(X_t,Y_t)\d t,\\
\d Y_t= \big(Z_t^{(2)}(X_t,Y_t)+ b_t(Y_t)\big)\d t +\si_t(Y_t)\d W_t,\ \ t\in [0,T],\end{cases}\end{equation}
where $[0,T]:=[0,\infty)$ when $T=\infty,$ $(W_t)_{t\in [0,T]}$ is an $m$-dimensional Brownian motion on a complete filtrated probability space
$(\OO, \F, \{\F_t\}_{t\in [0,T]},\P)$, and 
\beg{align*} &Z^{(i)}: [0,T]\times \R^{d_1+d_2}\to \R^{d_i},\ \ i=1,2,\\
&b: [0,T]\times\R^{d_2}\to \R^{d_2},\ \ \si: [0,T]\times \R^{d_2}\to \R^{d_2}\otimes\R^m\end{align*}
are measurable.  We assume that $Z_t^{(i)}(x,y)$ and $\si_t(y)$ are continuous in $(x,y)\in \R^{d_1+d_2}$ as $``$regular" coefficients, but  $b_t(y)$ only satisfies a local integrability condition in $(t,y)$ and is regarded as  a $``$singular" term. See assumptions $(A_1)$-$(A_3)$ below for details. 

 SDE \eqref{E1} is known as stochastic Hamiltonian system when the drifts are given by the gradients
of a Hamiltonian functional. The associated Fokker-Planck equation is called the Langevin equation or kinetic Fokker-Planck equation. In particular, when $Z_t^{(1)}(x,y)= y$, $X_t$ and $Y_t$ stand  for the location and speed of a fluid flow at time $t$ respectively.

In this paper, we investigate the well-posedness and exponential ergodicity of \eqref{E1} with drift $Z_t^{(2)}(x,y)+b_t(y)$ discontinuous in $(t,y)$. These two properties have been intensively
studied when the coefficients are regular enough or the invariant probability measure is known, which we summary as follows.

When $Z_t^{(1)}(x,y)$ and $Z_t^{(2)}(x,y)+b_t(y)$ satisfy a Dini-H\"older continuity condition in $(x,y)$, the well-posedness of \eqref{E1} has been proved in \cite{C16,WZ16},    see also the recent paper \cite{ZZZ} for the weak well-posedness of \eqref{E1} for  $d_1=d_2, Z^{(1)}(x,y)=y, \si=I_{d_2} ($the $d_2\times d_2$ identity matrix) and $Z^{(2)}+b$ being in   a weighted anisotropic Besov space.

When the SDE is time independent, the exponential ergodicity
for  special versions of \eqref{E1} has been studied in many references.  When the unique invariant probability measure is  given, the  hypocoercivity  introduced by Villani \cite{V} has attracted a lot of attentions ,
and has been further developed in a series of papers such as \cite{CHSG,  GS, GW19} based on an abstract analytic framework  built up  by Dolbeaut, Mouhot and Schmeiser \cite{DMS}. However, the study in this direction heavily relies on
the explicit formulation of the invariant probability measure, for which the drifts are given by (weighted) gradient of a Hamiltonian functional.  A crucial motivation in the ergodicity theory is to simulate the invariant probability measure by using the stochastic system. In this spirit, we study  the ergodicity  for the above general model with unknown invariant probability measure.  On the other hand, when the coefficients are regular enough,
the exponential ergodicity follows from some  modified dissipativity conditions, see for instance
\cite{W17} for the hypercontractivity implying the exponential ergodicity in related entropy, and \cite{RW21} for
an extension to McKean-Vlasov SDEs. 

In Section 2 and Section 3,   we investigate the well-posedness and exponential ergodicity for \eqref{E1} with $Z_t^{(2)}(x,y)$ only satisfying a local integrability condition 
in $(t,y)$, such that the corresponding results derived in the recent papers  \cite{R21, W21a,XZ} for non-degenerate SDEs are extended to the present degenerate setting. 
Finally, in Section 4 we extend the main results to McKean-Vlasov SDEs.

 \section{Well-posedness}

In this part we let $T\in (0,\infty)$ be finite. For any $r>0$ and $(x,y)\in \R^{d_1+d_2}$, let
$$ B_r(x,y):=\big\{(x',y')\in\R^{d_1+d_2}: |x-x'|+|y-y'|\le r\big\},\ \ B_r(y):=\big\{y'\in \R^{d_2}: |y-y'|\le r\big\}.$$ For any
$$(p,q)\in \scr K:=\Big\{(p,q):\ p,q\in (2,\infty), \ff {d_2} p+\ff 2 q<1\Big\},$$
define the $\tt L_q^p$-norm of a (vector or real valued) measurable function $f$ on $[0,T]\times\R^{d_2}$:
$$\|f\|_{\tt L_q^p}:= \sup_{y\in \R^{d_2}} \bigg(\int_0^T \|1_{B_1(y)}f_t\|_{L^p(\R^{d_2})}^q\d t\bigg)^{\ff 1 q}.$$
We write $f\in \tt L_q^p$ if $\|f\|_{\tt L_q^p}<\infty \mbox{ and } \tt L^p $ represents the $\tt L_q^p$-norm independent of $q.$ We use $\nn^{(1)}$ and $\nn^{(2)}$ to denote the gradient operators in $x\in \R^{d_1}$ and $y\in \R^{d_2}$ respectively, so that $(\nn^{(i)})^2$ is the corresponding Hessian operator. In case only one variable is concerned, we simply denote the gradient by $\nn$.

\beg{enumerate} \item[$(A_1)$] For any $n\ge 1$ there exists a constant $0<K_n<\infty$ such that
$$\max_{i=1,2}|Z_t^{(i)}(x,y)- Z_t^{(i)}(x',y')|\le K_n   (|x-x'|+|y-y'|),\ \ t\in [0,T], (x,y),(x',y')\in B_n(0).$$ Moreover,
$$\sup_{t\in [0,T], i=1,2} |Z_t^{(i)}(0)| <\infty. $$
\item[$(A_2)$]    $\si\si^*$ is invertible with $\|\si\|_\infty+\|(\si\si^*)^{-1}\|_\infty<\infty,$ and
$$\lim_{\vv\downarrow 0} \sup_{t\in [0,T], |y-y'|\le \vv} \|\si_t(y)-\si_t(y')\|=0.$$
Moreover, there exist $l\in\mathbb N$, $\{(p_i,q_i)\}_{0\le i\le l} $  and $1\le f_i\in \tt L_{q_i}^{p_i}, 0\le i\le l$ such that
$$|b|\le f_0,\ \ \|\nn \si\|\le \sum_{i=1}^l f_i.$$
\item[$(A_3)$]   There exist $\vv\in (0,1)$ and  $1\le V\in C^2(\R^{d_1+d_2})$ with
\beq\label{LY} \lim_{|x|+|y|\to\infty} V(x,y)=\infty,\ \ \limsup_{|x|+|y|\to\infty}\sup_{y'\in B_\vv(y)}
 \ff{|\nn^{(2)} V(x,y')|+ \|(\nn^{(2)})^2 V(x,y')\|}{V(x,y)}<\infty,\end{equation} such that
\beg{align*}
&\vv  \sup_{y'\in B_\vv(y)} \big\{|Z_t^{(1)}(x,y)|\|\nn^{(1)}  \nn^{(2)}  V(x,y')\|+ |Z_t^{(2)}(x,y)|\big(\|\nn^{(2)}V(x,y')\|+\| (\nn^{(2)})^2  V(x, y')\|\big)\big\} \\
&+\<Z_t^{(1)}(x,y), \nn V(\cdot,y)(x)\>+\<Z_t^{(2)}(x,y),\nn V(x,\cdot)(y)\>   \le \eta_t V(x,y)\end{align*}
holds for some $0\le \eta\in L^1([0,T]),$ all $t\in [0,T]$ and $(x,y)\in \R^{d_1+d_2}.$ \end{enumerate}

\beg{thm}\label{T1} Assume $(A_1)$ and $(A_2)$. Then for any initial value $(X_0,Y_0)$ the SDE $\eqref{E1}$ has a unique strong solution up to the life time $\zeta$. If $(A_3)$ holds then $\zeta=T$ and
\beq\label{EST} \E\Big[\sup_{t\in [0,T]} V(X_t,Y_t)\Big|\F_0\Big]\le c V(X_0,Y_0)\end{equation}
holds for some constant $0<c<\infty$.\end{thm}

By a truncation argument, we first consider the case where $\{K_n\}_{n\ge 1}$ is bounded. In this case, we may
follow the line of \cite{XXZZ} to prove the well-posedness by using Zvonkin's transform. The only difference is that in the present degenerate setting we apply this transform for the following time-dependent elliptic operators on $\R^{d_2}$:
$$L_t:=\ff 1 2 {\rm tr}\big(\si_t\si_t^*\nn^2\big)+b_t\cdot\nn,\ \ t\in [0,T].$$
Moreover, since $Y_t$ also depends on $X_t$, we have to reprove the Khasminskii estimate for $Y_t$ which is crucial in the proof of pathwise uniqueness.

\beg{lem}\label{L1} Assume $(A_2)$ and that  $\|Z^{(2)}\|_\infty <\infty.$
Then for any $(p,q)\in\scr K$ there exists an increasing function $H: [0,\infty)\to [0,\infty)$ such that
for any strong solution $(X_t,Y_t)_{t\in [0,T]}$ of $\eqref{E1}$,
$$\E\big[\e^{\int_0^T |f_t(Y_t)|^2\d t}\big]\le H(\|f\|_{\tt L_{q}^{p}}),\ \ f\in \tt L_{q}^{p}.$$
\end{lem}

\beg{proof} Let
$$\gg_s:= \big(\si_s^*(\si_s\si_s^*)^{-1}Z_s^{(2)}\big)(X_s,Y_s),\ \ s\in [0,T].$$
By $(A_2)$ and the boundedness of $Z^{(2)}$, we have
\beq\label{BP0} K:=\int_0^T \big\|\gg_s\big\|_\infty^2\d s<\infty.\end{equation}
So, by Girsanov's theorem,
$$\tt W_t:= W_t+\int_0^t \gg_s  \d s,\ \ t\in [0,T]$$
is an $m$-dimensional Brownian motion under the probability measure $\Q:=R\P$, where
$$R:=\e^{-\int_0^T\<\gg_s,\d W_s\>-\ff 1 2 \int_0^T|\gg_s|^2\d s}.$$
So, $Y_t$ solves the SDE
$$\d Y_t = b_t(Y_t)\d t +\si_t(Y_t)\d \tt W_t,\ \ t\in [0,T].$$
By \cite[Lemma 4.1]{XXZZ}, $(A_2)$ implies that
$$\E_\Q\e^{\int_0^T |f_t(Y_t)|^2\d t }\le H_0(\|f\|_{\tt L_{q}^{p}})$$
holds for some increasing function $H_0: [0,\infty)\to [0,\infty)$. Thus,
 \beg{align*} &\E\e^{\int_0^T |f_t(Y_t)|^2\d t }= \E_\Q\big[R^{-1} \e^{\int_0^T |f_t(Y_t)|^2\d t }\big]\\
&\le \big(\E[R^{-2}]\big)^{\ff 1 2}\Big(\E_\Q\e^{2\int_0^T |f_t(Y_t)|^2\d t }\big] \Big)^{\ff 1 2}
= \big(\E[R^{-2}]\big)^{\ff 1 2}\ss{H_0(2\|f\|_{\tt L_{q}^{p}})}.\end{align*}
Then the proof is finished since \eqref{BP0} implies
\beg{align*} &\E[R^{-2}] = \E \big[\e^{-2\int_0^T\<\gg_s,\d W_s\>+ \int_0^T|\gg_s|^2\d s}\big]\\
&\le \E \big[\e^{-2\int_0^T\<\gg_s,\d W_s\>- 2\int_0^T|\gg_s|^2\d s+3K}\big]=\e^{3K}<\infty.\end{align*}
 \end{proof}

 \beg{proof}[Proof of Theorem \ref{T1}]
 (a) The well-posedness up to life time. By a truncation argument,  in stead of $(A_1)$ we may and do assume that $Z_t^{(i)}$ are bounded and Lipschitz continuous in $(x,y)$ uniformly in $t\in [0,T],$ so that the life time of $(X_t,Y_t)$ is $T$.
  By \cite[Themore 3.2]{XXZZ},
 $(A_2)$ implies that for any $\ll\ge 0$ the PDE
\beq\label{PDE} (\pp_t+L_t)u_t=-\ll u_t-b_t,\ \ t\in [0,T], u_T=0\end{equation} for $u: [0,T]\times\R^{d_2}\to \R^{d_2}$  has a unique solution with
 $$\|u\|_\infty+\|\nn u\|_\infty +\|\nn^2 u\|_{\tt L_{q_0}^{p_0}}<\infty,$$ where $\nn$ and $\nn^2$ are the gradient and Hessian operators on $\R^{d_2}$. Moreover,   for   $\ee\in (0,1)$ there exits $\ll>0$ is large enough such that,
 \beq\label{W1} \|u\|_\infty+\|\nn u\|_\infty<\vv.\end{equation}
We now make  the Zvonkin's transform for $Y_t$:
$$\tt Y_t:= \Theta_t(Y_t),\ \ \Theta_t(y):=y+ u_t(y),\ \ t\in [0,T], y\in\R^{d_2}.$$ By \eqref{PDE} and generalised It\^o's formula (see \cite{XXZZ}), \eqref{E1} becomes
\beq\label{E1'}\beg{cases} \d X_t= \tt Z_t^{(1)} (X_t,\tt Y_t) \d t,\\
\d \tt Y_t= \tt Z_t^{(2)}(X_t,\tt Y_t) \d t +\tt \si_t(\tt Y_t)\d W_t,\ \ t\in [0,T],\end{cases}\end{equation}
where
\beg{align*} &\tt Z_t^{(1)}(x,y):= Z_t^{(1)}\big(x,\Theta_t^{-1}(y)\big),\\
&\tt Z_t^{(2)}(x,y):=\big((\nn \Theta_t)\circ\Theta_t^{-1}(y)\big) Z_t^{(2)}\big(x,\Theta_t^{-1}(y)\big)- \ll u_t\circ\Theta_t^{-1}(y),\\
&\tt\si_t(y):= \big((\nn \Theta_t)\circ\Theta_t^{-1}(y)\big)\si_t\circ\Theta_t^{-1}(y),\ \ t\in [0,T], (x,y)\in \R^{d_1+d_2}.\end{align*} Since $\|\nn u\|_\infty<1$, $\Theta_t$ is diffeomorphism so that the well-posedness of
\eqref{E1} is equivalent to that of \eqref{E1'}. Noting that the coefficients of \eqref{E1'} are bounded and continuous in $(x,y)$,  this SDE has a weak solution. By the Yamada-Watanable principle, it remains to prove the pathwise uniqueness of \eqref{E1'}.
This can be done as in \cite{XXZZ} by using Khasminskii's estimate in Lemma \ref{L1}. Below we present a detailed proof for completeness.

For any nonnegative measurable function $f$ on $\R^{d_2},$ consider its maximal functional
$$\scr M f(x):= \sup_{r\in (0,1)} \ff 1 {|B_r(0)|}\int_{B_r(0)} f(x+y)\d y,\ \ x\in \R^{d_2}.$$
By \cite[Lemma 2.1]{XXZZ}, there exists a constant $0<c<\infty$ such that for any function $f\in L^{\infty}(\R^d)$ with $\nn f\in L^1_{loc}(\R^d)$,
\beq\label{W2} \beg{split} &|f(x)-f(y)|\le c|x-y|\big(\|f\|_\infty+\scr M |\nn f|(x)+\scr M |\nn f |(y)\big),\ \ x,y\in\R^{d_2},\\
&\|\scr M |f|\|_{\tt L_q^p}\le c \|f\|_{\tt L_q^p}.\end{split}\end{equation}
Now, let $(X_t^{(i)},\tt Y_t^{(i)})_{t\in [0,T]}, i=1,2,$ be two solutions of \eqref{E1'} with
$(X_0^{(1)},\tt Y_0^{(1)})=(X_0^{(2)},\tt Y_0^{(2)}).$ By the Lipschitz continuity of $Z_t^{(i)}$ uniformly in $t$, \eqref{W2} and $\|\nn^2 u\|_{\tt L_q^p}<\infty$, we find   functions $1\le g_j\in \tt L_{q_j}^{p_j}\ (0\le j\le l)$ such that
\beg{align*} &\sum_{i=1,2}|\tt Z_t^{(i)}(x,y)-\tt Z_t^{(i)}(x',y')|^2 +\|\tt\si_t(y)-\tt\si_t(y')\|^2\\
&\le \big(|x-x'|^2+|y-y'|^2\big)\sum_{j=0}^l \big(g_j(t,y)^2+g_j(t,y')^2\big),\ \ t\in [0,T], (x,y),(x',y')\in \R^{d_1+d_2}.\end{align*}
Then by It\^o's formula, we find a constant $0<c<\infty$ and a martingale $M_t$ such that
$$\xi_t:= |X_t^{(1)}-X_t^{(2)}|^2+|\tt Y_t^{(1)}-\tt Y_t^{(2)}|^2,\ \ t\in [0,T]$$ satisfies
\beq\label{W3} \d \xi_t \le c \xi_t\sum_{j=0}^l \sum_{i=1,2} |g_j(t, \tt Y_t^{(i)})|^2+\d M_t,\ \ t\in [0,T], \xi_0=0.\end{equation}
By Lemma \ref{L1}, for all $0<\theta<\infty,$ we have
$$\E\big[\e^{\theta \int_0^{T} |g_j(t, \tt Y_t^{(i)})|^2\d t}\big]<\infty,\ \   0\le j\le l,~\\i=1,2.$$
So, by the stochastic Gronwall inequality, see \cite[Lemma 3.7]{XZ}, \eqref{W3} implies   $(X_t^{(1)}, \tt Y_t^{(1)})= (X_t^{(2)}, \tt Y_t^{(2)})$ for all $t\in [0,T].$

(b) Now, let $(A_3)$ hold, we aim to prove \eqref{EST}. By \eqref{W1}, we have $|Y_t-\tt Y_t|<\vv$.
Combining this with \eqref{LY} it suffices to prove \eqref{EST} for $\tt Y_t$ replacing $Y_t$, i.e.
\beq\label{EST'} \E\Big[\sup_{t\in [0,T]} V(X_t,\tt Y_t)\Big|\F_0\Big]\le c V(X_0,\tt Y_0)\end{equation}
holds for some constant $0<c<\infty$. By \eqref{W1}, \eqref{E1'}, and It\^o's formula, the boundedness of $\tt \si$ and $(A_3)$ imply that for some constant $0<C<\infty$,
\beq\label{E0*} \d V(X_t,\tt Y_t)\le C(1+\eta_t) V(X_t,\tt Y_t)\d t+\d M_t,\ \ t\in [0,T]\end{equation}  holds for some martingale $M_t$
with
$$\d \<M\>_t\le C V(X_t,\tt Y_t)^2\d t,\ \ t\in [0,T].$$
Let $\tau_n:=T\wedge \inf\{t\geq 0, |X_t^x|+|\tt Y_t^y|\geq n\}$ and the life time $$\zeta:=\lim_{n\to\infty}\tau_n,\ \ \ n\ge 1.$$
Then by the Gronwall Inequality we have
$$\E \big(V(X_{\tau_n}, \tt Y_{\tau_n})1_{\tau_n< T}\big|\F_0\big) \leq   V(X_0, \tt Y_0)\e^{c\int_0^T(1+\eta_s)\d s},$$
which implies
$$\E \big(1_{\tau_n< T}\big|\F_0\big) \leq \ff { V(X_0, \tt Y_0)\e^{c\int_0^T (1+\eta_s)\d s}} { \inf_{|x|+|y|\geq n}V(x,y)}.$$
Then by Fatou's Lemma, we have
$$\E \big(1_{\zeta\leq T}\big|\F_0\big) =0,$$
which further implies 
$$\P(\zeta\leq T)=\E \big(1_{\zeta\leq T}\big) =\E \Big(\E \big(1_{\zeta\leq T}\big|\F_0\big)\Big)=0.$$
With the definition of  the life time $\zeta, $ when $n\to \infty,$ we have $\zeta= T.$ 
Finally,  by a standard argument using the Burkholder-Davis-Gundy and Gronwall inequalities, we prove \eqref{EST'} for some constant $c>0$.
\end{proof}

 \section{Exponential ergodicity}

In this part we consider the time-homogeneous case such that \eqref{E1} becomes
\beq\label{E3}\beg{cases} \d X_t= Z^{(1)}(X_t,Y_t)\d t,\\
\d Y_t= \big\{Z^{(2)}(X_t,Y_t)+ b(Y_t)\big\}\d t +\si_t(Y_t)\d W_t,\ \ t\ge 0,\end{cases}\end{equation}
where
\beg{align*} &Z^{(i)}:   \R^{d_1+d_2}\to \R^{d_i},\ \ i=1,2,\ \
b:  \R^{d_2}\to \R^{d_2},\ \ \si:   \R^{d_2}\to \R^{d_2}\otimes\R^m\end{align*}
are measurable.  We investigate the ergodicity of the associated Markov process. To this end, we make the following assumption, which, according to Theorem \ref{T1}, implies the well-posedness and non-explosion of
this SDE.

\beg{enumerate} \item[$(B_1)$] For any $1\leq n <\infty$ there exists a constant $0<K_n<\infty$ such that
$$\sup_{i=1,2}|Z^{(i)}(x,y)- Z^{(i)}(x',y')|\le K_n   (|x-x'|+|y-y'|),\ \  (x,y),(x',y')\in B_n(0).$$ 
\item[$(B_2)$]    $\si\si^*$ is invertible with $\|\si\|_\infty+\|(\si\si^*)^{-1}\|_\infty<\infty,$ and there exists $p>(2\lor d)$ such that $|b|+\|\nn\si\|\in \tt L^p.$ 
\item[$(B_3)$]   There exist constants $\vv\in (0,1), 0<K<\infty$, an increasing function
$\Phi: [1,\infty)\to (0,\infty)$ with $\Phi(n)\to\infty$ as $n\to\infty$, and  $1\le V\in C^2(\R^{d_1+d_2})$ with
\beq\label{LY1} \lim_{|x|+|y|\to\infty} V(x,y)=\infty,\ \ \limsup_{|x|+|y|\to\infty}\sup_{y'\in B_\vv(y)}
 \ff{|\nn V(x,\cdot)(y')|+ \|\nn^2 V(x,\cdot)(y')\|}{\big(V \wedge \Phi(V)\big)(x,y)}=0,\end{equation}
 such that
\beg{align*}&\vv  \sup_{y'\in B_\vv(y)} \big\{|Z_t^{(1)}(x,y)|\|\nn^{(1)}  \nn^{(2)}  V(x,y')\|+ |Z_t^{(2)}(x,y)|(|\nn^{(2)} V(x,y')|+\| \nn^{(2)} \nn^{(2)}  V(x, y')\|)\big\} \\
&+\<Z_t^{(1)}(x,y), \nn V(\cdot,y)(x)\>+\<Z_t^{(2)}(x,y),\nn V(x,\cdot)(y)\>   \le K-\Phi(V(x,y)).\end{align*}
\end{enumerate}

By Theorem \ref{T1}, under $(B_1)$-$(B_3)$ the SDE \eqref{E3} is well-posed. Let $\{P_t\}_{t\ge 0}$ be the associated Markov semigroup, i.e.
$$P_tf(x,y)= \E[f(X_t^{x,y}, Y_t^{x,y})],\ \ f\in \B_b(\R^{d_1+d_2}), t\ge 0, (x,y)\in R^{d_1+d_2},$$ where $(X_t^{x,y}, Y_t^{x,y})$ solves 
\eqref{E3}  with initial value $(x,y)$. 
We investigate the ergodicity of $P_t$, i.e. it has a unique invariant probability measure $\mu$ such that
$$\lim_{t\to\infty } P_t^*\nu=\mu,\ \ \nu\in \scr P,$$ where $\scr P$ is the space of all  probability measures  on $\R^{d_1+d_2}$, and
$$(P_t^*\nu)(f):= \nu(P_t f)=\int_{\R^{d_1+d_2}} P_t f \d\nu,\ \ f\in \B_b(\R^{d_1+d_2}), \nu\in \scr P, t\ge 0.$$
In terms of $(B_3)$, we consider the convergence under the $V$-variation norm
$$\|\mu_1-\mu_2\|_{V}:=\sup_{f\in\B_b(\R^{d_1+d_2}), |f|\le V} |\mu_1(f)-\mu_2(f)|.$$
Under this norm, the space
$$\scr P_V:= \big\{\mu\in \scr P:\ \mu(V)<\infty\big\}$$ is a complete metric sapce. 
When $V\equiv 1$, we denote the norm  by $\|\cdot\|_{var}$ which is known as the total variation norm.
The Lyapunov condition $(B_3)$  implies the existence of invariant probability measure. 

\subsection{Main results and example}

To prove the ergodicity, we need the following assumption that any compact set is a petite  set of $P_t$:

\beg{enumerate} \item[$(B_4)$] Any compact set $D$ of $\R^{d_1+d_2}$ is $P_t$-petite, i.e. there exists $t_0>0$ and a non-trivial finite measure $\nu$ such that
$$\inf_{x\in D} P_{t_0}(x,\cdot)\ge \nu,$$
where $P_t(x,\cdot)$ is the transition probability kernel of $P_t$ at $x\in \R^{d_1+d_2}$ .
\end{enumerate}

\beg{thm}\label{T1.6.1} Assume $(B_1)$-$(B_3)$. Then the following assertions hold.
\beg{enumerate} \item[$(1)$] $P_t$ has an invariant probability measure $\mu$ such that
$$\mu\big(\Phi(\vv_0V)\big)=\int_{\R^{d_1+d_2}}\Phi\big(\vv_0V(x,y)\big)\mu(\d x,\d y)<\infty$$ holds for some $\vv_0>0$.
\item[$(2)$] If $P_t$ is $t_0$-regular for some $t_0>0$,  i.e. $\{P_{t_0}(x,\cdot): x\in\R^{d_1+d_2}\}$ are mutually equivalent, then
\beq\label{EGD} \lim_{t\to\infty} P_tf(x,y) =\mu(f),\ \ \mu\in \scr P,~f\in \B(\R^{d_1+d_2}). \end{equation}
 \item[$(3)$] If $(B_4)$ holds and $\Phi(r)\ge \dd r$ for some constant $\dd>0$ and all $r\ge 0$, then there exist constants $1<c<\infty,  \ll>0$
such that
 \beq\label{EX1}  \|P_t^*\mu_1-P_t^*\mu_2\|_{V}  \le c\e^{-\ll t} \|\mu_1-\mu_2\|_V,\ \ \mu_1,\mu_2\in \scr P_V, t\ge 0.\end{equation}
 Consequently, $\mu\in \scr P_V$ is the unique invariant probability measure of $P_t$, and 
$$  \|P_t^*\nu-\mu \|_{V}  \le c\e^{-\ll t} \|\nu-\mu\|_V,\ \ \nu \in \scr P_V, t\ge 0. $$
 \item[$(4)$] Let $(B_4)$ hold and $ H(r):=\int_0^r   \ff {\d s} {\Phi(s)} <\infty$ for $r\ge 0$.  If $\Phi$ is convex, then there exist constants $1<k<\infty,\ll>0$ such that
 \beq\label{EX0}  \|P_{t}^*\dd_{(x,y)}-\mu\|_{V}\le k \big\{1+H^{-1} (H(V(x,y))- k^{-1} t)\big\} \e^{-\ll t},\ \ (x,y)\in \R^{d_1+d_2}, t\ge 0,\end{equation}
 where $H^{-1}$ is the inverse of $H$ with $H^{-1}(r):=0$ for $r\le 0$.
  Consequently, if $H(\infty)<\infty$ then there exist constants $0<c,\ll,t^*<\infty$ such that
  \beq\label{EX2}\|P_t^*\mu_1-\mu_2\|_{V}\le c  \e^{-\ll t}\|\mu_1-\mu_2\|_{var},\ \ \forall t\ge t^*,  ~\mu_1,\mu_2\in   \scr P.\end{equation} \end{enumerate}
\end{thm}

In general, $(B_4)$ follows from a H\"ormander condition. In this spirit, we use the following explicit condition replacing $(B_4)$.

\beg{enumerate} \item[$(B_4')$] $d_1=d_2=d$, $\nn^{(2)}Z^{(1)}$ is invertible with 
$$\|\nn^{(2)}Z^{(1)}\|_\infty +\|(\nn^{(2)}Z^{(1)})^{-1}\|_\infty<\infty,$$ $\nn^{(2)} Z^{(1)}$ is H\"older continuous, and 
$\|\nn^{(1)} Z^{(1)}\|+\|(\nn^{(2)})^2 Z^{(1)}\|$ is locally bounded.
\end{enumerate}

\beg{thm}\label{TN} Assume $(B_1)$-$(B_3)$ and $(B_4')$, then   $(B_4)$ holds  and $P_{t}$ is $t_0$-regular for any $t_0>0$, so that all assertions of Theorem \ref{T1.6.1} apply.
\end{thm}

\paragraph{Example 3.1.} Simply consider  $d_1=d_2=d$ and that $\si_t(y)=I_{d\times d}$ is the identity matrix. We make the following choices of $b, Z^{(1)}$ and $Z^{(2)}$: 
\beg{enumerate} \item[$\bullet$] $b$ satisfies   $\|b\|_{\tt L^p}<\infty$ for some $p>d$.  For example,  it is easy to see that this   is true  when 
$$b(x):=\int_{\R^d} \frac{x-y}{|x-y|^{\alpha+1}}\nu(\mbox{d} y),\ \ \ x\in\R^d$$
for some $\aa\in (0,1)$ and a finite measure $\nu$ on $\R^d.$ This type drifts are of interests  in statistical physics, see \cite{12} and references therein. Here, we extend the existing study to degenerate setting. 
\item[$\bullet$] For some constants $c_1,c_2, c_3,\dd$ with $0<c_1, |c_2|, c_3<\infty$ and $0\leq\dd< \infty$, 
$$Z^{(1)}(x,y):= -c_1 (1+|x|)^\dd x + c_2 y,\ \ Z^{(2)}(x,y):= Z(x,y)-c_3 (1+|y|)^\dd y,$$
where  $Z: \R^{2d}\to\R^d$ is a locally Lispchitz continuous   with
$$\lim_{|(x,y)|\to\infty} \ff {|Z(x,y)|}{|(x,y)|}=0. $$
\end{enumerate} 
Take, for some constant $\theta\in (0,\infty),$  
$$V(x,y)= (1+|x|^2+|y|^2)^\theta,\ \ \ x,y\in\R^d.  $$  
Then $(B_1)$, $(B_2)$ and $(B_4')$ hold, so that    by  Theorem \ref{TN}, we have the following assertions. 
\beg{enumerate}
\item[$(1)$] When $\dd=0$ and $|c_2|$ is small enough, we find a constant $0<c_0<\infty$ such that $(B_3)$ holds for $\Phi(r)= c_0 r$.   
Assertions (1)-(3) in Theorem \ref{T1.6.1} imply that $P_t$ has a unique invariant probability measure  $\mu$ such that
 $\mu(|\cdot|^2)<\infty$, \eqref{EGD}  and \eqref{EX1} for some constants $0<c,\ll<\infty$ hold.
\item[$(2)$] When $\dd>0$, then $(B_3)$ holds for $\Phi(r)= c_0 (1+ r^{1+\dd/(2\theta)})$ for some constant $0<c_0<\infty$, so that   
  Theorem \ref{T1.6.1} (4) implies  \eqref{EX2} for some constants $c,\ll>0$. 
 \end{enumerate}

\subsection{Proofs} 

\beg{proof}
[Proof of Theorem \ref{T1.6.1}]   Once  (1) is proved,   (2) follows from  Doob's Theorem   \cite{Doob}. So, below we only prove (1), (3) and (4). 

(a) We  prove the existence of   invariant probability measure by using Zvonkin's transform.
Let 
$$L^0=\ff 1 2 {\rm tr} \Big( \si \si ^*\nn^2\Big)+ b\cdot\nn.$$
According to \cite[Lemma 2.5]{W21c}, $(B_2)$ implies that  there exists   $\lambda>0$ such that the PDE 
\beq\label{THT0} (L^0-\lambda)u=-b\end{equation} 
has a unique solution satisfying 
\beq\label{THT'} \|u\|_{\infty}+\|\nn u\|_{\infty}< \ee,\ \ \|\nn^2\|_{\tt L^p}<\infty,\end{equation}
where $\ee$ defined in $(B_2).$\\
Then  $\Theta(y):=y+u(y),~y\in\R^{d_2},$ gives rise to a diffeomorphism on $\R^{d_2}$. So, 
 \beq\label{THT}  \tt \Theta(x,y)=(x,\Theta(y)),\ \ (x,y)\in \R^{d_1+d_2}\end{equation}  is a  diffeomorphism on $\R^{d_1+d_2}.$

Let  $(X_t^{x,y},  Y_t^{x,y})_{t\geq 0}$ solve \eqref{E3} with initial value $(x,y)\in \R^{d_1+d_2}$, and let 
$$\hat{P}_t f(x,y):= \E [(X_t^{x,y},\Theta(Y_t^{x,y}))],\ \ t\ge 0, f\in \B_b(\R^{d_1+d_2}), (x,y)\in \R^{d_1+d_2}.$$ 
  Then $(P_t)_{t\geq 0}$ and $(\hat{P}_tf)_{t\geq 0}$ satisfy 
\beq\label{PHT} (\hat{P}_t f)(x,y)=\big( P_t(f\circ \tilde{\Theta})\big)(\tilde{\Theta}^{-1}(x,y)),~~(x,y)\in \R^{d_1+d_2}.\end{equation} 
So, $\mu$ is an invariant probability measure of $P_t$ if and only if 
\beq\label{pm}\hat \mu:= \mu\circ\tilde{\Theta}^{-1}\end{equation}
is  $\hat P_t$-invariant.   
Therefore,    it is sufficient to show that $ \hat{P}_t $ has an invariant probability measure. 
By  the  Bogoliov-Krylov theorem,  we only need to verify the tightness of  
\beq\label{pm1}\hat{\mu}_n:=\ff 1 n\int_0^n \hat{P}_s((0 ,0), \cdot)\d s,\ \ n\ge 1.\end{equation}
By $(B_3)$ and \eqref{THT'}, we find a constant $K\in (0,\infty)$ such that 
$$|V(x,\Theta(y))- V(x,y)|\le K \big\{V(x,\Theta(y))\land V(x,y)\big\},\ \ (x,y)\in \R^{d_1+d_2},$$
so that for $r_0:= \ff 1{1+K}\in (0,1)$, 
 \beq\label{W2'} \gamma_0 (V\circ \tilde{\Theta})(x, y)\leq V(x,y)\leq \gamma_0^{-1}(V\circ \tilde{\Theta})(x, y). \end{equation}
Combining these with $(B_3)$, \eqref{THT0} and applying It\^o's formula,  we find a constant $0<c_1<\infty$ such that
\beq\label{W3'}
\d (V\circ \tilde{\Theta})(X_t,  Y_t)\leq \Big(K-c_1 \Phi\big(\gamma_0(V\circ \tilde{\Theta})(X_t, Y_t)\big)\Big)\d t+\d M_t \end{equation}
for some martingale $(M_t)_{t\geq 0}$.  Letting $(X_0,Y_0)=(0,0),$ we deduce from \eqref{W3'} that 
\beg{align*}& \int_{\R^{d_1+d_2}}  \Phi\big(\gamma_0(V\circ \tilde{\Theta})  \big)  \d\hat \mu_n= \ff 1 n \int_0^n \E\Big[ \Phi\big(\gamma_0(V\circ \tilde{\Theta})(X_t^{(0,0)}, Y_t^{(0,0)})\big)\Big]\d t\\
&\le 
\ff {K+\ V\circ\tt\Theta(0)/n}{c_1} <\infty,\ \ n\ge 1.\end{align*} 
Since $ \Phi\big(\gamma_0(V\circ \tilde{\Theta})  \big)$ has compact level sets, this implies the tightness of $\{\hat\mu_n\}_{n\ge 1}$, and the weak limit $\hat\mu$ of a convergent subsequence gives an invariant  probability measure of $\hat P_t$. 
Moreover, 
$$\int_{\R^{d_1+d_2}}  \Phi\big(\gamma_0(V\circ \tilde{\Theta})  \big)  \d\hat \mu \le \ff{K}{c_1}<\infty.$$ 
Therefore, by \eqref{pm}  and \eqref{THT'}, $\mu:=\hat\mu\circ\tt\Theta$ is an invariant probability measure of $P_t$ and $\mu(\Phi(\vv_0V))<\infty$ holds for some constant $0<\vv_0<\infty$.

(b) In the situations of  (3) and (4), we have $\Phi(r)\ge c_0 r$ for some constant $0<c_0<\infty$ and all $r\ge 1$, so that (1) implies that $\mu(V)<\infty$. 
By   \eqref{PHT} and   the definition of weighted total variation norm, we obtain   
\beg{align*}&\|P_t^*\mu_1-P_t^*\mu_2\|_{V}=\|\hat{P}_t^*(\mu_1\circ \tilde{\Theta}^{-1})-\hat{P}_t^*(\mu_2\circ \tilde{\Theta}^{-1})\|_{V\circ \tilde{\Theta}^{-1}},\\
&|P_t^*\mu_1-\mu\|_{V}=\|\hat{P}_t^*(\mu_1\circ \tilde{\Theta}^{-1})-\hat{\mu}\|_{V\circ \tilde{\Theta}^{-1}},\ \ \mu_1,\mu_2\in \scr P_V,\end{align*}  where $\mu$ and $\hat{\mu}$ are the above constructed invariant probability measures of $(P_t)_{t\geq 0}$ and $(\hat{P}_t)_{t\geq 0}$, respectively.  Combining this with \eqref{W2'}, we derive 
\beg{align*} 
&\|P_t^*\mu_1-P_t^*\mu_2\|_{V}\leq \gamma_0^{-1}\|\hat{P}_t^*(\mu_1\circ \tilde{\Theta}^{-1})-\hat{P}_t^*(\mu_2\circ \tilde{\Theta}^{-1})\|_V, \\
&\|P_t^*\mu_1-\mu\|_{V}\leq \gamma_0^{-1}\|\hat{P}_t^*(\mu_1\circ \tilde{\Theta}^{-1})-\hat{\mu}\|_{V},\ \ \mu_1, \mu_2\in \scr P_V.\end{align*}
So, it remains to verify   \eqref{EX1} and \eqref{EX0}  for $\hat{P}^*_t$ replacing $P^*_t$. 

By \eqref{W3'} and $\Phi(r)\geq c_0r,$ we have
$$\d V(x,\Theta(y))\leq \big(K-c_1c_0\gamma_0V(x,\Theta(y))\big)\d t+\d M_t.$$
By invoking $(B_4)$ and \eqref{PHT},     any compact subset of $\R^{d_1+d_2}$  is $\hat P_t$-petite.  Furthermore,  from \eqref{W3'} we deduce that $(\hat{P}_t^*)_{t\geq 0}$ admits the Lyapunov condition: for some constant $0<k_1,~k_2=c_1c_0\gamma_0<\infty$ and for any $(x,y)\in \R^{d_1+d_2}$,
$$\hat{P}_tV(x,\Theta(y))=\E\big(V(x,\Theta(y))\big)\leq \frac{k_1}{k_2}+\e^{-k_2t}V(x,\Theta(y)).$$
Consequently,  the Harris theorem \cite[Theorem 1.3]{HM} and \eqref{W2'} yield that there exist constants $0<c,\lambda<\infty$ such that 
$$\|\hat{P}^*_t\delta_{(x,y)}-\hat{\mu}\|_V\leq \e^{-\lambda t}\|\delta_{(x,y)}-\hat{\mu}\|_V\leq c\e^{-\lambda t}V(x,y),~~t\geq 0. $$
Thus,  by following the line of part (c) in the proof of  \cite[Theorem 2.1]{W21c}, we find constants $c,\ll>0$ such that 
\beq\label{W6}\|\hat{P}^*_t\mu_1-\hat{P}^*_t\mu_2\|_V\leq c\e^{-\lambda t}\|\mu_1-\mu_2\|_V,\ \ \mu_1,\mu_2\in \scr P_V. \end{equation}
This immediately implies \eqref{EX1} for $\hat{P}^*_t$ replacing $P^*_t$. 

Next,  by  \eqref{W6}, the semigroup property of $(\hat{P}_t)_{t\geq 0}$ and the invariance of $\hat{\mu}$, we have
$$\|\hat{P}^*_t\delta_{(x,y)}-\hat{\mu}\|_V=\|\hat{P}^*_{t/2}\hat{P}^*_{t/2}\delta_{(x,y)}-\hat{P}^*_{t/2}\hat{\mu}\|_V\leq c\e^{-\lambda t/2}\|\hat{P}^*_{t/2}\delta_{(x,y)}-\hat{\mu}\|_V.$$
On the other hand,  by the proof of \cite[(2.35)]{W21c}, $(B_4)$ implies 
$$\hat{P}_{2/t}V(x,y)\leq c_2 \Big( 1+H^{-1}(H(V(x,y))-\frac{t}{2c_2}\Big),~~0<c_2<\infty$$
so that 
$$\|\hat{P}^*_t\delta_{(x,y)}-\hat{\mu}\|_V\leq c\e^{\frac{-\lambda t}{2}}\|\hat{P}^*_{2/t}\delta_{(x,y)}-\hat{\mu}\|_V \leq c_3 \e^{\frac{-\lambda t}{2}}\big(\hat{P}^*_{2/t}\delta_{(x,y)}(V)+\hat{\mu}(V)\big)$$
$$ \leq c_4 \Big( 1+H^{-1}(H(V(x,y))-\frac{t}{2 c_2}\Big)\e^{-\lambda t},  0<c_3,c_4<\infty.$$
Therefore,    \eqref{EX0} holds for $\hat{P}^*_t$ replacing $P^*_t$.  
\end{proof}

\beg{proof}[Proof of Theorem \ref{TN}] By Proposition \ref{P1} below,   $(B_1)$-$(B_3)$   and $(B_4')$ imply $(B_4)$ for any $t_0>0$ and 
$\nu(\d x, \d y):= \inf_{(x',y'), (x'',y'')\in D} p_{t_0}(x,y; x',y') 1_{D}(x,y)\d x\d y$. .\end{proof} 

In the following Proposition \ref{P1},    $(B_3)$ is weakened as  

\beg{enumerate} 
\item[$(B_3')$]   There exist constants $\vv\in (0,1), K>0$, and  $1\le V\in C^2(\R^{d_1+d_2})$ with
\beq\label{LY1} \lim_{|x|+|y|\to\infty} V(x,y)=\infty,\ \ \limsup_{|x|+|y|\to\infty}\sup_{y'\in B_\vv(y)}
 \ff{|\nn V(x,\cdot)(y')|+ \|(\nn^{2})^2 V(x, y')\|}{V(x,y)}<\infty,\end{equation}
 such that
\beg{align*}&\vv  \sup_{y'\in B_\vv(y)} \big(|Z^{(1)}(x,y)|\|\nn^{(1)}  \nn^{(2)}  V(x,y')\|+ |Z^{(2)}(x,y)|(|\nn^{(2)} V(x,y')|+\| (\nn^{2})^2 V(x, y')\|)\big) \\
&+\<Z^{(1)}(x,y), \nn V(\cdot,y)(x)\>+\<Z^{(2)}(x,y),\nn V(x,\cdot)(y)\>   \le K V(x,y).\end{align*}
\end{enumerate}

\beg{prp}\label{P1} Assume that $(B_1), (B_2), (B_3')$   and $(B_4')$ hold.  If $\si$  is H\"older continuous, then $P_t$ has a heat kernel $p_t(x,y;x',y')$ with respect  to the Lebesgue measure such that
$$\inf_{(x,y), (x',y')\in B_k(0)} p_t(x,y; x',y')>0,\ \forall t,k>0.$$\end{prp}

 To prove this result, we apply the Harnack inequality presented in \cite{GIMV} for the PDE
 \beq\label{PDE''} \beg{split} \pp_t f_t (x,y)= &\,- y\cdot\nn^{(1)} f_t(x,y) - \big(Z^{(2)} \cdot\nn^{(2)} f_t\big)(x,y) \\
 &+  {\rm div}^{(2)}\big(a\nn^{(2)} f_t\big)(x,y)+ (Uf_t)(x,y),\ \ t\ge 0, x,y\in \R^{d},\end{split}\end{equation}
 where ${\rm div}^{(2)}$ is the divergence operator in the second component $y$, and 
 $$Z^{(2)}: \R^{2d} \to \R^d,\ \ U: \R^{2d}\to\R,\ \ a: \R^{2d} \to \R^{d\otimes d}$$ are measurable satisfying the following assumption.
 
 \beg{enumerate} 
\item[$(B_2')$]    $a$ is invertible, 
  $\nn^{(2)} a$ exists,  such that $ |Z^{(2)}|+|U|+\|a\|+\|a^{-1} \|$ is locally bounded in $\R^{2d}.$ \end{enumerate} 
 
 The following Harnack inequality is essentially due to \cite{GIMV}.
 
 \beg{lem}\label{LN1} Assume $(B_2')$. Then there exists a constant $r_0\in (0,1)$ such that for any $t>0$ and $r\in (0,r_0]$, there exists a locally bounded function
 $$\varphi_{t,r}: \R^{2d}\times \R^{2d}\to (0,\infty)$$
 such that any  positive weak solution (in the sense of integration by parts) $f_t$ of \eqref{PDE''} satisfies 
 $$f_t(x',y')\le \varphi_{t,r}(x,y;x',y') f_{t+r} (x,y),\ \ (x,y), (x',y') \in \R^{2d}.$$
 \end{lem} 

\beg{proof} By \cite[Theorem 3 and Remark 4]{GIMV}, there exist small constants $r_0,r_1\in (0,1)$, such that for any $t>0$ there exists a locally bounded function
$$C_t: (0,r_0]\times \R^{2d}\to (0,\infty)$$
such that any positive solution $f$ of \eqref{PDE''} satisfies
$$f_t(x,y)\le C_{t,r}(x,y) \inf_{(x',y')\in B(x,y;r_1)} f_{t+r}(x', y'),\ \ \ r\in (0, r_0], $$
where $B(x,y; r_1):=\{(x',y')\in \R^{2d}:\ |(x-x',y-y')|\le r_1\}.$
For any $(x,y), (x',y')\in \R^{2d},$ let  
$$n=n(x,y;x',y'):=\inf\big\{n\in \mathbb N:\ |(x-x',y-y')|\le  nr_1\big\},$$ and denote
$$(x_i,y_i):= (x,y)+ \ff i n (x'-x, y'-y),\ \ 0\le i\le n.$$
Then 
$$f_t(x,y)\le C_{t,\ff{r} n}(x_0,y_0)f_{t+\ff r n}(x_1,y_1)\le\cdots\le f_{t+r} (x',y') \prod_{i=0}^{n-1} C_{t+ \ff {ir}n, \ff r n} (x_i,y_i).$$
Therefore, the desired estimate holds for
$$\varphi_{t,r}(x,y;x',y'):=  \prod_{i=0}^{n-1} C_{t+ \ff {ir}n, \ff r n} (x_i,y_i),\ \ n:=n(x,y;x',y').$$
 \end{proof} 
 
 Next, we extend Lemma \ref{LN1} to the following more general PDE:
  \beq\label{PDE'} \beg{split} \pp_t f_t (x,y)= &\,- \big(Z^{(1)}  \cdot\nn^{(1)} f_t+ Z^{(2)} \cdot\nn^{(2)} f_t\big)(x,y)  \\
  & +   {\rm div}^{(2)}\big(a \nn^{(2)} f_t\big)(x,y)+(Uf_t)(x,y),\ \ t\ge 0,\ x,y\in\R^d.\end{split} \end{equation}

 \beg{lem}\label{LN2} Assume $(B_2')$ and $(B_4')$. Then there exists a constant $r_0\in (0,1)$ such that for any $t>0$ and $r\in (0,r_0]$, there exists a locally bounded function
 $$\varphi_{t,r}: \R^{2d}\times \R^{2d}\to (0,\infty)$$
 such that any  positive weak solution (in the sense of integration by parts) $f_t$ of \eqref{PDE'} satisfies 
 $$f_t(x',y')\le \varphi_{t,r}(x,y;x',y') f_{t+r} (x,y),\ \ (x,y), (x',y') \in \R^{2d}.$$
\end{lem} 

\beg{proof} To transform \eqref{PDE'} into \eqref{PDE''}, we make the change of variable
$$(x,y)\mapsto (x,\tt y):= \big(x, (Z^{(1)}(x,y))^{-1} \big).$$
Let  $\phi(x,\cdot):= \big(Z^{(1)}(x,\cdot)\big)^{-1}$ and
$$\tt f_t(x,\tt y):= f_t(x,\phi(x,\tt y)),\ \  \ \tt U(x,\tt y):= U(x,\phi(x,\tt y)),\ \ t\ge 0, x,\tt y\in\R^d.$$
Then   \eqref{PDE'} implies
\beq\label{R1} \beg{split} \pp_t \tt f_t (x,\tt y)= &\,- \Big(Z^{(1)}  \cdot\nn^{(1)} f_t+ Z^{(2)} \cdot\nn^{(2)} f_t\Big)\big(x,\phi(x,\tt y)\big)  \\
  & +   {\rm div}^{(2)}\Big(a \nn^{(2)} f_t\Big)\big(x,\phi(x,\tt y)\big) +(Uf_t)\big(x,\phi(x,\tt y)\big),\ \ t\ge 0,\ x,y\in\R^d.\end{split} \end{equation}
By chain rule, we obtain
\beg{align*}&\Big( \nn^{(1)}\tt f_t \Big)(x,\tt y) = (\nn^{(1)} f_t)\big(x,\phi(x,\tt y)\big)+ \Big(\nn^{(1)}\phi(x,\tt y)\Big)\big(\nn^{(2)}f_t\big)\big(x,\phi(x,\tt y)\big),\\
&\Big(\nn^{(2)}\tt f_t\Big)(x,\tt y) = \Big(\nn^{(2)}\phi(x,\tt y)\Big)\big(\nn^{(2)} f_t\big)\big(x,\phi(x,\tt y)\big),\\
&\Big(\big(\nn^{(2)}\big)^2\tt f_t\Big)(x,\tt y)=\Big(\nn^{(2)}\phi(x,\tt y)\Big)^2 \Big(\big(\nn^{(2)}\big)^2 f_t \Big)\big(x,\phi(x,\tt y)\big) \\
&\qquad\qquad \qquad \qquad +\Big((\nn^{(2)})^2 \phi(x,\tt y)\Big) \big(\nn^{(2)}f_t\big)\big(x,\phi(x,\tt y)\big).\end{align*}
So,
\beq\label{R2} \beg{split} & \big(\nn^{(2)}f_t\big)\big(x,\phi(x,\tt y)\big) = \Big(\nn^{(2)}\phi(x,\tt y)\Big)^{-1} \big(\nn^{(2)}\tt f_t\big)(x,\tt y),\\
& \big(\nn^{(1)} f_t\big)\big(x,\phi(x,\tt y)\big) = \big(\nn^{(1)}\tt f_t\big)(x,\tt y) - \Big(\nn^{(1)} \phi(x,\tt y)\Big) \Big(\nn^{(2)} \phi(x,\tt y)\Big)^{-1} \big(\nn^{(2)} \tt f_t\big)(x,\tt y).\end{split}\end{equation} 
So, letting
$$A^{x,\tt y}:= \Big(\nn^{(2)} \phi(x,\tt y)\Big)^{-1},$$
we derive 
\beq\label{R3} \beg{split}  & {\rm div}^{(2)}\big(a \nn^{(2)}  f_t\big)\big(x,\phi(x,\tt y)\big)=\sum_{i,j=1}^d \big[\pp_i^{(2)} \big(a_{ij} \pp_j^{(2)} f_t\big)\big](x,\phi(x,\tt y))\\
&= \sum_{i,j,k=1}^d A^{x,\tt y}_{ik} \pp_k^{(2)} \big((a_{ij}\pp_j^{(2)}f_t)(x,\phi(x,\tt y))\big)\\
&= \sum_{i,j,k,l=1}^d A_{ik}^{x,\tt y} \pp_k^{(2)}\big(a_{ij}(x,\phi(x,\tt y)) A_{jl}^{x,\tt y} \pp_l^{(2)}\tt f_t (x,\tt y)\big)\\
&= \sum_{k,l=1}^d \pp_{k}^{(2)}\Big( \big[(A^{x,\tt y})^* a(x,\phi(x,\tt y)) A^{x,\tt y} \big]_{kl}\pp_l^{(2)}\tt f_t(x,\tt y)\Big) \\
&\quad -\sum_{i,j,k,l=1}^d \big(\pp_k^{(2)}A_{ik}^{x,\tt y}\big) a_{ij}(x,\phi(x,\tt y)) A_{jl}^{x,\tt y} \pp_l^{(2)}\tt f_t(x,\tt y)\\
&= {\rm div}^{(2)}\big(\tt a \nn^{(2)}\tt f_t\big)(x,\tt y) - \big(\hat Z\cdot \nn^{(2)}\tt f_t\big)(x,\tt y),\end{split}\end{equation}
where  
$$\tt a(x,\tt y):= \big(A^{x,\tt y}\big)^*a(x,\phi(x,\tt y)) A^{x,\tt y},\ \ \hat Z_l:=\sum_{i,j,k=1}^d \big(\pp_k^{(2)} A^{x,\tt y}_{ik} \big)a_{ij}(x,\phi(x,\tt y)) A_{jl}^{x,\tt y},\ 1\le l\le d.$$
Substituting \eqref{R2} and \eqref{R3} into \eqref{R1}, and noting that $Z_1(x,\phi(x,\tt y))=\tt y$, we derive
  \beg{align*} \pp_t \tt f_t (x,\tt y)= &\,- \big(\tt y \cdot\nn^{(1)} \tt f_t+ \tt Z^{(2)} \cdot\nn^{(2)} \tt f_t\big)(x,\tt y)  \\
  & +   {\rm div}^{(2)}\big(\tt a \nn^{(2)} f_t\big)(x,\tt y)+(\tt U\tt f_t)(x,y),\ \ t\ge 0,\ x,y\in\R^d,\end{align*}
  for the above defined $\tt U$, $\tt a$ and 
$$\tt Z^{(2)}(x,\tt y):=   \hat Z(x,\tt y)+\big(\nn^{(1)}\phi(x,\tt y)\big)^*\tt y+ \big(\nn^{(2)}\phi(x,\tt y)\big)^*Z^{(2)}(x,\phi(x,\tt y)),\ \ x,\tt y\in\R^d.$$
Combining this with $(B_2')$ and $(B_4')$,  we may apply Lemma \ref{LN1} to this PDE to derive the desired estimate.
\end{proof} 

We also need the following result for the existence of heat kernel. 

\beg{lem}\label{LN3} Assume that $(B_1), (B_2)$   and $(B_4')$ hold,  and the solution to $\eqref{E3}$ is non-explosive. Then \eqref{E3} has heat kernel $p_t $;   namely, for any $t>0$ and the solution $(X_t, Y_t)$ starting at $(x_0,y_0)$, the distribution of $(X_t,Y_t)$ has a density $p_t(x_0,y_0;\cdot)$ with respect to the Lebesgue measure.\end{lem}

\beg{proof} (a) We first assume $b=0$ and $Z^{(1)}(x,y)=y$, but allow $\si$ also depends on $x$ such that $\si$ is H\"older continuous and 
$\|\si^*\|_\infty+\|(\si\si^*)^{-1}\|_\infty<\infty$. When    $Z^{(2)}$ is bounded, 
the existence of heat kernel follows from \cite[Theorem 1.5]{LPP}. In general, for any $n\ge 1$ let
$$Z^{2,n}(x,y):= Z^{(2)}(\varphi_n(x), \varphi_n(y)),\ \ \varphi_n(x):= x1_{\{|x|\le n\}} + \ff {nx}{|x|} 1_{\{|x|>n\}}.$$
Let $(X_t^{(n)}, Y_t^{(n)})$ solve the SDE 
$$\beg{cases} \d X_t^{(n)}= Y_t^{(n)} d t,\\
\d Y_t^{(n)}= Z^{2,n}(X_t^{(n)}, Y_t^{(n)})\d t+\si(Y_t)\d W_t,\ \ (X_0,Y_0)=(x_0,y_0).\end{cases}$$
Then for any $t>0$ and $n\ge 1$, the distribution of $(X_t^{(n)}, Y_t^{(n)})$ is absolutely continuous with respect to the Lebesgue measure; i.e. for any null set $A\subset \R^{2d}$,
$\P((X_t^{(n)}, Y_t^{(n)})\in A)=0.$ Letting
$$\tau_n:=\inf\{t\ge 0: |X_t|\lor |Y_t|\ge n\},$$
we have $(X_t,Y_t)= (X_t^n, Y_t^n)$ for $t\le \tau_n$. By the non-explosion we have $\tau_n\to\infty$ as $n\to\infty$, so that
$$\P((X_t,Y_t)\in A)\le \lim_{n\to\infty}\big\{ \P((X_t^{(n)}, Y_t^{(n)})\in A) + \P(\tau_n<t)\big\}=0$$
holds for all null set $A$. Thus, the heat kernel exists. 

(b) Let $b=0$ and for $Z^{(1)}$ satisfying $(B_4')$. As shown in the proof of Lemma \ref{LN2}, with the transform $(x,y)\mapsto \big(x, (Z^{(1)}(x,y))^{-1}\big)$ we reduce the situation (a),  so that the heat kernel exists.
Finally, when $\|b\|_{\tt L^p}<\infty$ for some $p>d$, by \cite[Lemma 2.5]{W21c}, when $\ll>0$ is large enough, 
the  PDE 
$$\ff 12 {\rm tr} \big\{\si\si^*\nn^2\big\} u+b\cdot \nn u= -b+\ll u$$
for $u: \R^d\to \R^d$  has a unique solution such that 
\beq\label{ASS} \|\nn^2 u\|_{\tt L^p}<\infty,\ \ \|\nn u\|_\infty \le \ff 1 2.\end{equation} 
Moreover, by the Sobolev embedding theorem, $\|\nn^2 u\|_{\tt L^p}<\infty$  implies that $\nn u$ is H\"older continuous. 
By It\^o's formula, we see that  $(X_t,\tt Y_t)= (X_t, \Theta(Y_t))$ solves the SDE
$$\beg{cases} \d X_t= \tt Z^{(1)}(X_t,\tt Y_t)\d t,\\
\d \tt Y_t= \tt Z^{(2)}(X_t, \tt Y_t)\d t+\tt \si(Y_t)\d W_t,\ \ (X_0,\tt Y_0)=(x_0,\Theta(y_0)),\end{cases}$$
where 
\beg{align*}& \tt Z^{(1)}(x,\tt y):= Z^{(1)}(x, \Theta^{-1}(\tt y)),\\
&\tt Z^{(2)}(x,\tt y):=\big((\nn \Theta)  Z^{(2)}(x,\cdot)\big)(\Theta^{-1}(\tt y))+ \ll u(\Theta^{-1}(\tt y)),\\
&\tt\si(\tt y):=\big((\nn\Theta)\si\big)(\Theta^{-1}(\tt y)),\ \ x,\tt y\in \R^d.\end{align*}
Thus,   $(X_t, \tt Y_t):= (X_t, \Theta(Y_t))$ solves the SDE of type \eqref{E3} with $b=0$, so that by Step (a) it has a heat kernel.
By \eqref{ASS}, this implies that $(X_t,Y_t)$ has heat kernel as well.  
\end{proof} 


\beg{proof}[Proof of Proposition \ref{P1}]  
By Lemma \ref{LN3},   for any $t>0$, $(X_t,Y_t)$ has a distribution density function (heat kernel) $p_t(x_0,y_0;\cdot)$. 
Let $a=\ff 1 2 \si\si^*$, so that 
$$\ff 1 2{\rm tr} \big(\si\si^*\nn^{(2)}\big)f= {\rm div}^{(2)}\big(a\nn^{(2)} f\big) +\big( {\rm div}^{(2)}a\big)\cdot \nn^{(2)} f,~~f\in C^2$$
where $({\rm div }^{(2)}a)_i:= \sum_{j=1}^d \pp_j^{(2)} a_{ij}.$ 
As shown in Step (b) in the proof of Lemma \ref{LN3}, with Zvonkin's transform $(X_t, \tt Y_t):= (X_t, \theta(Y_t))$ for $b+  {\rm div}^{(2)}a$ replacing $b$, we may and do assume that
$b+  {\rm div}^{(2)}a=0$, so that  the generator of $(X_t,Y_t)$ becomes 
$$Lf:= {\rm div}^{(2)}\big(a \nn^{(2)}f\big)+ Z^{(1)}\cdot \nn^{(1)}f + Z^{(2)}\cdot \nn^{(2)}f.$$
It is easy to see that the adjoint operator of $L$ in $L^2(\R^{2d})$ is
$$L^*f=  {\rm div}^{(2)}\big(a \nn^{(2)}f\big)-Z^{(1)}\cdot \nn^{(1)}f -Z^{(2)}\cdot \nn^{(2)}f -\big({\rm div}^{(1)}Z^{(1)} +{\rm div}^{(2)}Z^{(2)}\big)f,$$
 Hence, that the heat kernel $f_t:=   p_t(x_0,y_0;\cdot)$ solves the equation \eqref{PDE'} for 
 $$U:= \big({\rm div}^{(1)}Z^{(1)} +{\rm div}^{(2)}Z^{(2)}\big).$$
 So,    by   Lemma \ref{LN2} with $\beta=\beta_0\land \ff t 2$, we derive 
$$p_t(x_0,y_0; x,y)\ge C_t(x,y; x',y') p_{t-\beta} (x_0,y_0; x', y'),\ \ (x,y), (x',y')\in \R^{2d}$$
for some function $C_t: \R^{2d}\times \R^{2d}\to (0,\infty)$ satisfying
$$c_t(\gamma):=\inf_{(x,y),(x',y')\in B_{\gamma}(0)} c_t(x,y; x', y')>0.$$ 
Let $k>0$ be a constant. For any $\gamma\ge k$, we obtain
\beg{align*} &\inf_{ (x,y)\in B_k(0)} p_t(x_0,y_0; x,y)\ge c_t(\gamma) \sup_{(x', y')\in B_{\gamma}(0)} p_{t-\beta} (x_0,y_0; x', y')\\
&\ge \tt c_t(\gamma)  \int_{B_N(0)} p_{t-\beta}(x_0,y_0; x' ,y')\d x' \d y'
= \tt c_t(\gamma) \P\big((X_t,Y_t)\in B_{\gamma}(0)\big),\end{align*} 
where $\tt c_t(\gamma):= \ff{c_t(N)}{|B_{\gamma}(0)|}>0$, for $|B_{\gamma}(0)|$ the volume of $B_{\gamma}(0)$. 
On the other hand, by Theorem \ref{T1}, $(B_1), (B_2)$ and $(B_3')$ imply \eqref{EST}, so that  for
$k_{\gamma}:=\inf_{(x,y)\notin B_{\gamma}(0)} U(x,y)$, 
$$ \P\big((X_t,Y_t)\notin B_{\gamma}(0)\big)\le \P\big(U(X_t,Y_t)\ge k_{\gamma}\big)\le \ff{c}{k_{\gamma}} U(x_0,y_0).$$
Taking large enough $N$ such that $ \ff{c}{k_{\gamma}} \sup_{(x_0,y_0)\in K} U(x_0,y_0)\le \ff 1 2$, we obtain 
$$\inf_{(x_0,y_0), (x,y)\in K} p_t(x_0,y_0; x,y)\ge \tt c_t(\gamma) \inf_{(x_0,y_0)\in K} \Big(1-  \ff{c}{k_{\gamma}} U(x_0,y_0)\Big)\ge \ff 1 2 \tt c_t(\gamma)>0.$$
So, the desired assertion holds.

\end{proof} 

 \section{Extension to McKean-Vlasov SDEs}

In this part we extend Theorem \ref{T1} and Theorem \ref{TN}  to McKean-Vlasov SDEs. Consider the following distribution dependent SDE on $\R^{d_1+d_2}$:
\begin{equation}\label{E41}
\beg{cases} \d X_t=  Z_t^{(1)}(X_t, Y_t)\d t,\\
\d  Y_t=  \big( Z_t^{(2)}(X_t,  Y_t,  \L_{(X_t,  Y_t)})+b_t(Y_t)\big)\d t+ \si(Y_t)\d W_t,\ \ t\in [0,T],\end{cases}\end{equation}
where $ \L_{(X_t,  Y_t)}\in\scr P,~t\in[0,T]$ is the law of $(X_t,  Y_t),$ $Z^{(1)}, b, \si$ and $W$ are as in \eqref{E1}, and
$$Z^{(2)}: [0,T]\times\R^{d_1+d_2}\times \scr P\to \R^{d_2}$$ is measurable. 
 We first study the well-posedness of  \eqref{E41}, then investigate the uniform ergodicity for the time-homogeneous model
\begin{equation}\label{EQ4.1}
\beg{cases} \d X_t=  Z^{(1)}(X_t, Y_t)\d t,\\
\d  Y_t=  \big( Z^{(2)}(X_t,  Y_t,  \L_{(X_t,  Y_t)})+b(Y_t)\big)\d t+ \si(Y_t)\d W_t,\ \ t\ge 0.\end{cases}\end{equation}

\subsection{Well-posedness of \eqref{E41} } 

Let $\hat{\scr P}$ be a sub-space of $\scr P$. We call \eqref{E41} well-posed for distributions in $\hat{\scr P}$, if for any $\F_0$-measurable random variable $(X_0,Y_0)$ with $\L_{(X_0,Y_0)}\in \hat{\scr P}$ (respectively, any initial distribution $\gg\in \hat{\scr P}$), \eqref{E41} has a unique strong solution (respectively, unique weak solution) $(X_t,Y_t)$ such that
$$[0,T]\ni t \to \L_{(X_t,Y_t)}\in \hat{\scr P}$$ is continuous in the weak topology. 

To extend Theorem \ref{T1}, let $\dd_0$ be the Dirac measure at $0\in \R^{d_1+d_2}$ and denote 
$$Z_t^{(2)} (x,y):= Z_t^{(2)}(x,y,\dd_0),\ \ t\in [0,T], (x,y)\in \R^{d_1+d_2}.$$
Let  
$$  C_b^w([0,T];\hat{\scr P}):=\big\{\mu\in[0,T]\rightarrow \hat{\scr P} \mbox{ is weakly continous,}  \sup_{t\in [0,T]} \mu_t(V)<\infty\big\}. $$

\beg{enumerate} \item[$(\hat A)$] $(A_1)$-$(A_3)$ hold for the above defined $Z^{(2)}_t(x,y)$. Moreover,   for any $n\ge 1$ and any   $\mu\in C_b^w([0,T];\hat{\scr P}), $  there exists a constant $K_{n,\mu}>0$ such that  
$$|Z_t^{(2)}(x,y,\mu_t)-Z_t^{(2)}(x',y', \mu_t)|\le K_{n,\mu} |(x-x',y-y')|,\ \ (x,y), (x',y')\in B_n(0).$$  
\end{enumerate} 
The following result extends Theorem \ref{T1} to the distribution dependent setting as well as \cite[Theorem 1.1]{R21} to the present degenerate case. 

\beg{thm}\label{T4.1} Assume $(\hat A)$ for $\hat {\scr P}= \scr P$ or  $\hat {\scr P}= \scr P_V$.  \beg{enumerate} 
\item[$(1)$] If $\hat {\scr P}=\scr P$, and there exists $0\le K\in L^2([0,T])$ such that
\beq\label{CD} \big|Z_t^{(2)}(x,y, \mu)- Z_t^{(2)}(x,y, \nu)\big|\le K_t \|\mu-\nu\|_{var},\ \  \mu,\nu\in \scr P,\end{equation}
then  $\eqref{E41}$ is well-posed for distributions in $\scr P$.
\item[$(2)$] Let $\hat {\scr P}=\scr P_V.$ If    
\beq\label{CD'} \beg{split} & |Z_t^{(2)}(x,y, \mu)- Z_t^{(2)}(x,y, \nu)|\sup_{y'\in B_\vv(y)}  \big(1+  |\nn^{(2)} V(x,y')| +\|(\nn^{(2)})^2 V(x,y')\|\big) \\
&\le K_t \|\mu-\nu\|_{V},\ \  (x,y)\in \R^{d_1+d_2},\  \mu,\nu\in \scr P_V,\end{split}\end{equation}
then  $\eqref{E41}$ is well-posed for distributions in $\scr P_V$.\end{enumerate}
\end{thm}

\beg{proof} Let $(X_0,Y_0)$ be $\F_0$-measurable with $\gg:= \L_{(X_0,Y_0)}\in \hat{\scr P}$. For any 
$$\mu\in \C^\gg:=\{\mu\in C_b^w([0,T]; \hat{\scr P}): \mu_0=\gg\},$$ 
 assumption   $(\hat A)$ together with  \eqref{CD} or \eqref{CD'}  implies $(A_1)$-$(A_3)$ for $Z_t^{(2)}(x,y,\mu_t)$ replacing $Z_t^{(2)}$, 
 so that by Theorem \ref{T1},  the SDE
 \beq\label{DE} \beg{cases} \d X_t^\mu=  Z_t^{(1)}(X_t^\mu, Y_t^\mu)\d t,\\
\d  Y_t^\mu=  \{ Z_t^{(2)}(X_t^\mu,  Y_t^\mu,  \mu_t)+b_t(Y_t^\mu)\}\d t+ \si(Y_t^\mu)\d W_t,\ \ t\in [0,T],\end{cases}\end{equation} is well-posed, 
where $ (X_0^\mu,Y_0^\mu)=(X_0,Y_0).$ 
By  \cite[Theorem 3.1]{HRW}, it suffices to show that 
$$\Psi:  \C^\gg\to \C^\gg,\ \ \Psi_t(\mu):=  \L_{(X_t^\mu, Y_t^\mu)}$$ 
has a unique fixed point in $\scr C^{\gamma}$. Below we prove this for  $\hat {\scr P}= \scr P$ and  $\hat {\scr P}= \scr P_V$ respectively. 

(1) Let $\hat {\scr P}=\scr P$ and \eqref{CD} holds. For  $\mu,\nu\in \C^\gg$, we reformulate \eqref{DE} as 
\beq\label{DE'} \beg{cases} \d X_t^\mu=  Z_t^{(1)}(X_t^\mu, Y_t^\mu)\d t,\\
\d  Y_t^\mu=  \big( Z_t^{(2)}(X_t^\mu,  Y_t^\mu,  \nu_t)+b_t(Y_t^\mu)\big)\d t+ \si(Y_t^\mu)\d \tt W_t,\ \ t\in [0,T], \end{cases}\end{equation}
where  $(X_0^\mu,Y_0^\mu)=(X_0,Y_0)$ and 
\beg{align*} &\tt W_t:= W_t-\int_0^t     \xi_s^{\mu,\nu}\d s,\\
&\xi_s^{\mu,\nu}:=  \big(\si_s^* (\si_s\si_s^*)^{-1}\big)(Y_s^{\mu}) \big(Z_s^{(2)}(X_s^\mu,  Y_s^\mu,  \nu_s) - Z_s^{(2)}(X_s^\mu,  Y_s^\mu,  \mu_s)\big).\end{align*}
By $(A_2)$ and \eqref{CD}, there exists a constant $0<c_1<\infty$ such that
\beq\label{XIS} |\xi^{\mu,\nu}_s|^2\le c_1 K_s^2\|\mu_s-\nu_s\|_{var}^2,\ \ s\in [0,T].\end{equation}
Since $[0,T]\ni s\mapsto \|\mu_s-\nu_s\|_{var}$ is measurable and bounded,  by  Girsanov's theorem,  $\tt W$ is a Brownian motion under the weighted probability measure $\Q:= R_T\P$, where
$$R_t:= \e^{\int_0^t\<\xi^{\mu,\nu}_s,\d W_s\>-\ff 1 2 \int_0^t |\xi^{\mu,\nu}_s|^2\d s},\ \ \ t\in [0,T]$$ is a martingale. Then by the weak uniqueness of  \eqref{DE}, the law of $(X_t^\mu, Y_t^\mu)$ 
under $\Q$ satisfies 
\beq\label{CD2} \L_{(X_t^\mu,Y_t^\mu)|\Q} =\L_{(X_t^\nu,Y_t^\nu)}= \Psi_t(\nu),\ \ t\in [0,T].\end{equation} 
Combining this with the martingale property of $R_t$,   Pinsker's inequality, \eqref{XIS} and letting $\E_\Q$ be the expectation with respect to $\Q$, we derive
\beq\beg{split} \label{CDN} &\|\Psi_t(\mu)-\Psi_t(\nu)\|_{var}^2 := \sup_{|f|\le 1} \big| \E\big[ f(X_t^\mu,Y_t^\mu) (1-R_t)\big]\big|^2 \le \big(\E[|1-R_t|]\big)^2\\
&\le  2\E[R_t\log R_t]= 2\E_\Q [\log R_t]= \int_0^t \E_\Q[|\xi^{\mu,\nu}_s|^2]\d s \le c_1 \int_0^t K_s^2 \|\mu_s-\nu_s\|_{var}^2\d s.\end{split}\end{equation}
This implies that when $\ll>0$ is large enough, $\Psi$ is contractive under the complete metric
$$\rr_\ll(\mu,\nu):=\sup_{t\in [0,T]} \e^{-\ll t}\|\mu_t-\nu_t\|_{var},\ \ \mu,\nu\in \C^\gg,$$
and hence has a unique fixed point. 

(2) Let $\hat {\scr P}=\scr P_V$ and \eqref{CD} holds. The following argument is similar to the proof of \cite[Theorem 1.1 (1)]{R21}, we include here for completeness.
 Let 
 $$\C_N^\gg:=\Big \{\mu\in \C^\gg: \sup_{t\in [0,T]} \mu_t(V) \e^{-Nt} \le N(1+\gg(V))\Big\},\ \ N\ge 1.$$ It suffices to find a constant $N\ge 1$ such that   $\Psi \C_N^\gg\subset \C_N^\gg$ and $\Psi$ has a unique fixed point in $\C_N^\gg.$

Firstly,   by $(A_3)$ for $Z_t^{(2)}(x,y):= Z_t^{(2)}(x,y,\dd_0)$ and \eqref{CD'}, we find a constant $0<c_0<\infty$  such that 
\beg{align*}
&\vv  \sup_{y'\in B_\vv(y)} \Big(|Z_t^{(1)}(x,y)|\|\nn^{(1)}  \nn^{(2)}  V(x,y')\|+ |Z_t^{(2)}(x,y,\mu_t)|\big(\|\nn^{(2)}V(x,y')\|+\| (\nn^{(2)})^2  V(x, y')\|\big)\Big) \\
&+\<Z_t^{(1)}(x,y), \nn V(\cdot,y)(x)\>+\<Z_t^{(2)}(x,y,\mu_t),\nn V(x,\cdot)(y)\>   \le  \eta_t  V(x,y)+c_0K_t \mu_t(V) \end{align*}
holds for  all $t\in [0,T]$,  $(x,y)\in \R^{d_1+d_2}$ and $ \mu\in \C^\gg$. 
As in \eqref{E0*},  by combining  this with  $(A_2)$ and It\^o's formula for $$(X_t^\mu, \tt Y_t^\mu):= (X_t^\mu,\Theta_t(Y_t^\mu)),\ \ t\in [0,T],$$ 
  we find a constant $0<c_1<\infty$ such that
\beq\label{CD0} \d  \big(V(X_t^\mu,\tt Y_t^\mu)\big)^2  \le c_1\Big((1+ \eta_t)  \big(V(X_t^\mu,\tt Y_t^\mu)^2+  K_t^2 \mu_t(V)^2\big)  \Big)  \d t + \d M_t,\ \ t\in [0,T],\end{equation} 
 for some martingale $M_t$. Then there exists a constant $0<c_2<\infty$ such that 
 \beq\label{CDS}\beg{split} & \E[V(X_t^\mu, \tt Y_t^\mu)^2|\F_0 ] \le \e^{c_1\int_0^t (1+ \eta_r)\d r} V(X_0,Y_0)^2  + \int_0^t c_1  K_s^2 \e^{c_1\int_s^t (1+  \eta_r)\d r} \mu_s(V)^2 \d s \\
& \le c_2  V(X_0,Y_0)^2 + c_2  N^2\big(1+\gg(V) \big)^2 \int_0^t K_s^2 \e^{2Ns}\d s\\
 &\le c_2 V(X_0,Y_0)^2 +c_2\big( N\big(1+\gg(V) \big) \e^{Nt} \big)^2    \int_0^t K_s^2 \e^{-2N(t-s)}\d s,\ \ t\in [0,T], \ \ \mu\in \C_N^\gg.\end{split}\end{equation}
 Noting that  $K\in L^2([0,T])$ implies
 $$\lim_{N\to\infty}\sup_{t\in [0,T]} \int_0^t K_s^2\e^{-2N(t-s)}\d s =0,$$ 
 while   $(A_3)$, \eqref{W1} and Jensen's inequality  yield 
 $$\big(\Psi_t(\mu)\big)(V):= \E [V(X_t^\mu,Y_t^\mu)]\le c_3 \E[V(X_t^\mu,\tt Y_t^\mu)]\le c_3 \E\Big[\big(\E(V(X_t^\mu, \tt Y_t^\mu)^2|\F_0)\big)^{\ff 1 2} \Big] $$ holds for some constant $0<c_3<\infty$, 
 we find   a constant $1\leq N_0<\infty$ such that \eqref{CDS} yields 
\beg{align*} &\sup_{t\in [0,T]} \{\Psi_t(\mu)\}(V)\e^{- Nt} \le  c_2 \E [V(X_0,Y_0) ]+ c_2   N \big(1+\gg(V)\big)   \int_0^t K_s^2 \e^{-2N(t-s)}\d s\\
 &= c_2 \gg(V)  + c_2   N \big(1+\gg(V)\big)   \int_0^t K_s^2 \e^{-2N(t-s)}\d  \le N(1+\gg(V)),\ \ N\ge N_0.\end{align*} 
 Therefore, $\Psi \C_N^\gg \subset \C_N^\gg$ for $N\ge N_0.$ 

 Next,  let $N\ge N_0$. We intend to prove that $\Psi$ has a unique fixed point in $\C_N^\gg$, by  using the Girsanov  transform defined  in (1) to show that for large $\ll>0$,  $\Psi$ is contractive in
 the following  complete metric on $\C_N^\gg$:
 $$\rr_{V,\ll}(\mu,\nu):= \sup_{t\in [0,T]} \e^{-\ll t} \|\mu_t-\nu_t\|_V.$$
  By \eqref{CD'} and $(A_2)$, 
 we find a constant $c_0 >0$ such that instead of \eqref{XIS}, 
\beq\label{XIS'} |\xi_s^{\mu,\nu}|^2 \le c_0 K_s^2   \|\mu_s-\nu_s\|_V^2,\ \ s\in [0,T],\ \mu,\nu\in \C_N^\gg.\end{equation} 
 So,    this together with   \eqref{CD2}, \eqref{CDS} and  \eqref{EST} yields that for some constants $0<c_1(N), c_2(N) <\infty$, 
\beq\label{RPP2} \beg{split} &\|\Phi_t(\mu)-\Phi_t(\nu)\|_V = \sup_{|f|\le V} \big|\E\big[f(X_t^{\mu},Y_t^{\mu})(1-R_t)\big]\big|\le \E\big[V(X_t^{\mu},Y_t^{\mu})|1-R_t| \big]\\
&\le \E\Big[\big(\E(V(X_t^\mu,Y_t^\mu)^2|\F_0 )\big)^{\ff 1 2} \big(\E [|R_t-1|^2|\F_0 ]\big)^{\ff 1 2}\Big]\\
   &\le c_1(N) \E \Big[V(X_0,Y_0) \big(\E [R_t^2-1|\F_0]\big)^{\ff 1 2}\Big],\ \ \mu,\nu\in\C_N^\gg,\end{split} \end{equation} 
   and due to $\e^{r}-1\le r \e^r $  for   $r\in \R$, 
    \beg{align*} &\E [R^2_t-1|\F_0]\le \E \bigg[\e^{2 \int_0^t \<\xi_s^{\mu,\nu},\d W_s\>- \int_0^t |\xi_s^{\mu,\nu}|^2\d s} -1\bigg|\F_0\bigg]\\
&\le \E \bigg[\e^{2 \int_0^t \<\xi_s^{\mu,\nu},\d W_s\>- 4\int_0^t |\xi_s^{\mu,\nu}|^2\d s+ \int_0^t 3c_0K_s^2 \|\mu_s-\nu_s\|_V^2\d s} \bigg|\F_0\bigg] -1 \\
&= \e^{  \int_0^t 3c_0K_s^2 \|\mu_s-\nu_s\|_V^2\d s}   -1 \le \e^{  \int_0^t 3c_0K_s^2 \|\mu_s-\nu_s\|_V^2\d s}  \int_0^t  3c_0 K_s^2\|\mu_s-\nu_s\|_V^2\d s\\
&\le    c_2(N)   \int_0^t K_s^2 \|\mu_s-\nu_s\|_V^2\d s,\ \ \ \mu,\nu\in \C_N^\gg,\ t\in [0,T].\end{align*}   
Combining this with \eqref{RPP2}, we find a constant $0<c_3(N)<\infty$ such that 
\beg{align*} &\rr_\ll(\Phi(\mu), \Phi(\nu))=\sup_{t\in [0,T]} \e^{-\ll t} \|\Phi_t(\mu)-\Phi_t(\nu)\|_V\\
&\le c_3(N) (1+\gg(V)) \rr_\ll(\mu,\nu)\sup_{t\in [0,T]}\bigg(\int_0^t K_s^2\e^{-2\ll(t-s)}\d s\bigg)^{\ff 1 2},\ \ \mu,\nu\in \C_N^\gg, \ t\in [0,T].\end{align*} 
Since 
$$\lim_{\ll\to\infty} \sup_{t\in [0,T]}\bigg(\int_0^t K_s^2\e^{-2\ll(t-s)}\d s\bigg)^{\ff 1 2}=0,$$
 when $\ll>0$ is large enough, $\Psi$ is contractive on the complete metric space  $(\C_N^\gg, \rr_\ll),$ so that it has a unique fixed point in $\C_N^\gg$ as desired. 

\end{proof} 

\subsection{Uniform ergodicity of  \eqref{EQ4.1}}
Assume that
\beg{enumerate} \item[$(\hat B)$] $(B_1)$-$(B_3)$ and $(B'_4)$ hold  for $ Z^{(2)}(x,y):= Z^{(2)}(x,y,\dd_0)$. 
Moreover, there exists a constant $0<\kk<\infty$ such that for all $(x,y)\in\R^{d_1+d_2}$ and $\gamma_1,  \gamma_2\in \scr P,$
$$|Z^{(2)}(x,y,\gamma_1)-Z^{(2)}(x,y,\gamma_2)|\leq \kk \|\gamma_1-\gamma_2\|_{var}.$$
\end{enumerate}

To investigate the uniform exponential ergodicity, we consider the following reference SDE for $\gg\in \scr P,$
\begin{equation}\label{EQ4.2}
\beg{cases} \d X_t^{\gamma}=  Z^{(1)}(X_t^{\gamma}, Y_t^{\gamma})\d t,\\
\d  Y_t^{\gamma}=  \big( Z^{(2)}(X_t^{\gamma},  Y_t^{\gamma},  \gamma)+b(Y_t^{\gamma})\big)\d t+ \si(Y_t^{\gamma})\d W_t.\end{cases}\end{equation}

Let $(P_t^{\gamma})^*\nu=\L_{(X_t^{\gamma},Y_t^{\gamma})}$ with $\L_{(X_0^{\gamma},Y_0^{\gamma})}=\nu.$ Denote by $(X_t^{\gamma,x,y},Y_t^{\gamma,x,y})$ the solution to \eqref{EQ4.2} with the initial value $ (x,y)\in\R^{d_1+d_2}.$

\beg{thm}\label{Thm4} Assume $(\hat B)$. If $\kk$ is small enough and 
   $\Phi$ in $(B_3)$ is convex with $\int_0^{\infty} \frac {1}{ \Phi(s)}\d s < \infty,$ then   $P_t^*$ associated with \eqref{EQ4.1} has a unique invariant probability measure $\hat\mu\in \scr P$ such that $\hat\mu(\Phi(\ee_0V))<\infty$ for some $\ee_0>0$,  and there exists constants $c,\lambda>0$ such that 
   $$\|P_t^*\nu-\hat\mu\|_{var}\leq ce^{-\lambda t}\|\hat\mu-\nu\|_{var},\ \ t\ge 0, \nu\in \scr P.$$
\end{thm}

\begin{proof}
According to \cite[Lemma 3.3]{W21c}, it is sufficient to find   constants $0<k_0,c,\ll<\infty$ such that  when $\kk<k_0$,  for any $\mu \in \scr P, ~ (P_t^{\gamma})_{t\geq 0}$ has a unique invariant measure $\mu_{\gamma}$ satisfying 
\begin{equation}\label{EQ4.3}
\|(P_t^{\gamma})^*\mu-\mu_{\gamma}\|_{var}\leq ce^{-\lambda t}\|\mu-\mu_{\gamma}\|_{var}.
\end{equation}
By the uniform Harris type theorem \cite[Lemma 3.3]{W21c},  we only need find   $t_0, t_1>0$ and a measurable set $B\in \scr B(\R^{d_1+d_2}),$ such that 
\begin{equation}\label{EQ4.4}
\inf_{\gamma\in \scr P, z\in \R^{d_1+d_2} }P_{t_0}^{\gamma}(z,B)>0,
\end{equation}
\begin{equation}\label{EQ4.5}
\sup_{\gamma\in \scr P; z,z'\in B}\|(P_{t_1}^{\gamma})^* \dd_z  -(P_{t_1}^{\gamma})^*\dd_{z'} \|_{var}<2.
\end{equation}
Below we prove these two estimates for $t_1=1, B= B_k(0)$ for large enough $k>0$ and some $t_0>1$. 

(a) Proof of \eqref{EQ4.4}.  Let $Z_t^{(2)}(x,y)=Z_t^{(2)}(x,y,\dd_0)$,  $(X_t^{\dd_0, x,y}, Y_t^{\dd_0, x,y})_{t\geq 0}$ solve \eqref{E3} with  initial value $(x,y)\in \R^{d_1+d_2},$
and let $(P_t^{\dd_0})_{t\geq 0}$ be the associated Markov semigroup. 
By Theorem \ref{T1.6.1} (1) and (3), $P_t^{\dd_0}$ has a unique invariant probability measure $\mu$ such that  \eqref{EX2} holds for some constants $0<c,\ll<\infty$. Consequently,  
there exists a constant $1\leq t_0<\infty$ such that
$$\|(P_t^{\dd_0})^*\nu- \mu\|_{var}\le \ff 1 4,\ \ t\ge t_0,\ \nu\in \scr P.$$
Taking $k>0$ such that $\mu(B_k(0))>\ff 3 4$, this implies that for $B:=B_k(0)$, 
\beq\label{FY1}   P_{t_0}^{\dd_0}1_{B}(x,y)\ge 1- P_{t_0}^{\dd_0} 1_{B^c}(x,y)\ge 1- \mu(B^c)-\ff 1 4 \ge \ff 1 2,\ \ (x,y)\in \R^{d_1+d_2}.\end{equation}
Now, for any $\gg\in \scr P$, let 
$$\xi_s^\gg:= \big(\si^*(\si\si^*)^{-1}\big)(X_s^{\dd_0,x,y},Y_s^{\dd_0,x,y}) \big(Z_2^{(2)}(X_s^{\dd_0,x,y},Y_s^{\dd_0,x,y}, \gg)- Z_2^{(2)}(X_s^{\dd_0,x,y},Y_s^{\dd_0,x,y},\dd_0)\big).$$
Then $(B_1)$ and $(B_2)$ imply 
$$|\xi_s^\gg|\le c_1,\ \ s\ge 0,\ \gg\in\scr P $$
for some constant $0<c_1<\infty$. Let
$$R:= \e^{\int_0^{t_0} \<\xi_s^\gg,\d W_s\>-\ff 1 2 \int_0^{t_0}|\xi_s|^2\d s}.$$
Thus,
$$\E[R^{-1}]=\E[\e^{-\int_0^{t_0} \<\xi_s^\gg,\d W_s\>+\ff 1 2 \int_0^{t_0}|\xi_s|^2\d s}]\le \e^{c_1^2 \frac{t_0}{2}},$$
and by Girsanov's theorem,
\beq\label{FY0} \L_{(X_t^{\dd_0,x,y}, Y_t^{\dd_0,x,y})|\Q}= \L_{(X_t^{\gg,x,y}, Y_t^{\gg,x,y})},\ \ t\in [0,t_0].\end{equation} 
So that   Schwarz inequality and \eqref{FY1} yield 
$$P_{t_0}^\gg 1_{B}(x,y) = \E [1_{B}(X_{t_0}^{\dd_0,x,y}, Y_{t_0}^{\dd_0, x,y}) R] \ge \ff{(\E[1_{B}(X_{t_0}^{\dd_0,x,y}, Y_{t_0}^{\dd_0,x,y})])^2 }{\E [R^{-1}]}\ge \ff 1 2 \e^{-c_1^2t_0}>0,$$ for any 
$(x,y)\in\R^{d_1+d_2}$ and any $\gg\in \scr P.$ 
Hence, \eqref{EQ4.4} holds.

(b) Proof of \eqref{EQ4.5}: Recall that $B=B_k(0)$. Let $p_1(x,y;\cdot)$ be the distribution density of $P_1^{\dd_0}(x,y;\cdot).$ By Proposition \ref{P1},  we have
$$\dd:= \inf_{(x,y),(x',y')\in B} p_1(x,y;x',y') >0,$$
where $p_1$ is the heat kernel of $P_1$. Then for any $z,z'\in B$, 
\beq\label{FY2}\beg{split} & \|(P_1^{\dd_0})^* \dd_z-(P_1^{\dd_0})^*\dd_{z'}\|_{var}= \int_{\R^{2d}}|p_1(z;x,y)-p_1(z';x,y)|\d x\d y\\
&\le \int_{\R^{2d}}\big(p_1(z;x,y)+ p_1(z';x,y) - 2[p_1(z;x,y)\land  p_1(z';x,y)]\big)\d x\d y \\
&\le 2 -2 \dd {\rm vol}(B):=\dd'<2.\end{split}\end{equation} 
On the other hand, by $(\hat B)$, \eqref{FY0} and Pinsker's inequality as in \eqref{CDN}, we find a constant $0<c_2<\infty$ such that 
$$\|(P_1^{\dd_0})^* \dd_z- (P_1^\gg)^*\dd_z\|_{var}\le c_2\kk,\ \  \ z\in\R^{d_1+d_2}.$$
Combining this with  \eqref{FY2} we conclude that when $\kk$ is small enough, \eqref{EQ4.5} holds for $t_1=1. $
 
\end{proof}


\end{document}